\input amstex
\magnification\magstephalf
\documentstyle{amsppt}

\hsize 5.72 truein
\vsize 7.9 truein
\hoffset .39 truein
\voffset .26 truein
\mathsurround 1.67pt
\parindent 20pt
\normalbaselineskip 13.8truept
\normalbaselines
\binoppenalty 10000
\relpenalty 10000
\csname nologo\endcsname 


\font\bc=cmb10
\font\tenbsy=cmbsy10

\catcode`\@=11

\def\myitem#1.{\item"(#1)."\advance\leftskip10pt\ignorespaces}

\def\qedsymbol{{\mathsurround\z@$\square$}}
\redefine\qed{\relaxnext@\ifmmode\let\next\@qed\else
  {\unskip\nobreak\hfil\penalty50\hskip2em\null\nobreak\hfil
    \qedsymbol\parfillskip\z@\finalhyphendemerits0\par}\fi\next}
\def\@qed#1$${\belowdisplayskip\z@\belowdisplayshortskip\z@
  \postdisplaypenalty\@M\relax#1
  $$\par{\lineskip\z@\baselineskip\z@\vbox to\z@{\vss\noindent\qed}}}
\outer\redefine\beginsection#1#2\par{\par\penalty-250\bigskip\vskip\parskip
  \leftline{\tenbsy x\bf#1. #2}\nobreak\smallskip\noindent}
\outer\redefine\genbeginsect#1\par{\par\penalty-250\bigskip\vskip\parskip
  \leftline{\bf#1}\nobreak\smallskip\noindent}

\def\next{\let\@sptoken= }\def\next@{ }\expandafter\next\next@
\def\@futureletnext#1{\let\nextii@#1\futurelet\next\@flti}
\def\@flti{\ifx\next\@sptoken\let\next@\@fltii\else\let\next@\nextii@\fi\next@}
\expandafter\def\expandafter\@fltii\next@{\futurelet\next\@flti}

\let\zeroindent\z@
\let\savedef@\endproclaim\let\endproclaim\relax 
\define\chkproclaim@{\add@missing\endroster\add@missing\enddefinition
  \add@missing\endproclaim
  \envir@stack\endproclaim
  \edef\endit@{\leftskip\the\leftskip\rightskip\the\rightskip}}
\let\endproclaim\savedef@
\def\thing@{.\enspace\egroup\ignorespaces}
\def\thingi@(#1){ \rm(#1)\thing@}
\def\thingii@\cite#1{ \rm\@pcite{#1}\thing@}
\def\thingiii@{\ifx\next(\let\next\thingi@
  \else\ifx\next\cite\let\next\thingii@\else\let\next\thing@\fi\fi\next}
\def\thing#1#2#3{\chkproclaim@
  \ifvmode \medbreak \else \par\nobreak\smallskip \fi
  \noindent\advance\leftskip#1
  \hskip-#1#3\bgroup\bc#2\unskip\@futureletnext\thingiii@}
\let\savedef@\endproclaim\let\endproclaim\relax 
\def\endit{\endproclaim\endit@\let\endit@\undefined}
\let\endproclaim\savedef@
\def\defn#1{\thing\parindent{Definition #1}\rm}
\def\lemma#1{\thing\parindent{Lemma #1}\sl}
\def\prop#1{\thing\parindent{Proposition #1}\sl}
\def\thm#1{\thing\parindent{Theorem #1}\sl}
\def\cor#1{\thing\parindent{Corollary #1}\sl}
\def\conj#1{\thing\parindent{Conjecture #1}\sl}

\def\remk#1{\thing\zeroindent{Remark #1}\rm}

\def\narrowthing#1{\chkproclaim@\medbreak\narrower\noindent
  \it\def\next{#1}\def\next@{}\ifx\next\next@\ignorespaces
  \else\bgroup\bc#1\unskip\let\next\narrowthing@\fi\next}
\def\narrowthing@{\@futureletnext\thingiii@}

\def\@cite#1,#2\end@{{\rm([\bf#1\rm],#2)}}
\def\cite#1{\in@,{#1}\ifin@\def\next{\@cite#1\end@}\else
  \relaxnext@{\rm[\bf#1\rm]}\fi\next}
\def\@pcite#1{\in@,{#1}\ifin@\def\next{\@cite#1\end@}\else
  \relaxnext@{\rm([\bf#1\rm])}\fi\next}

\advance\minaw@ 1.2\ex@
\atdef@[#1]{\ampersand@\let\@hook0\let\@twohead0\brack@i#1,\z@,}
\def\brack@{\z@}
\let\@@hook\brack@
\let\@@twohead\brack@
\def\brack@i#1,{\def\next{#1}\ifx\next\brack@
  \let\next\brack@ii
  \else \expandafter\ifx\csname @@#1\endcsname\brack@
    \expandafter\let\csname @#1\endcsname1\let\next\brack@i
    \else \Err@{Unrecognized option in @[}%
  \fi\fi\next}
\def\brack@ii{\futurelet\next\brack@iii}
\def\brack@iii{\ifx\next>\let\next\brack@gtr
  \else\ifx\next<\let\next\brack@less
    \else\relaxnext@\Err@{Only < or > may be used here}
  \fi\fi\next}
\def\brack@gtr>#1>#2>{\setboxz@h{$\m@th\ssize\;{#1}\;\;$}%
 \setbox@ne\hbox{$\m@th\ssize\;{#2}\;\;$}\setbox\tw@\hbox{$\m@th#2$}%
 \ifCD@\global\bigaw@\minCDaw@\else\global\bigaw@\minaw@\fi
 \ifdim\wdz@>\bigaw@\global\bigaw@\wdz@\fi
 \ifdim\wd@ne>\bigaw@\global\bigaw@\wd@ne\fi
 \ifCD@\enskip\fi
 \mathrel{\mathop{\hbox to\bigaw@{$\ifx\@hook1\lhook\mathrel{\mkern-9mu}\fi
  \setboxz@h{$\displaystyle-\m@th$}\ht\z@\z@
  \displaystyle\m@th\copy\z@\mkern-6mu\cleaders
  \hbox{$\displaystyle\mkern-2mu\box\z@\mkern-2mu$}\hfill
  \mkern-6mu\mathord\ifx\@twohead1\twoheadrightarrow\else\rightarrow\fi$}}%
 \ifdim\wd\tw@>\z@\limits^{#1}_{#2}\else\limits^{#1}\fi}%
 \ifCD@\enskip\fi\ampersand@}
\def\brack@less<#1<#2<{\setboxz@h{$\m@th\ssize\;\;{#1}\;$}%
 \setbox@ne\hbox{$\m@th\ssize\;\;{#2}\;$}\setbox\tw@\hbox{$\m@th#2$}%
 \ifCD@\global\bigaw@\minCDaw@\else\global\bigaw@\minaw@\fi
 \ifdim\wdz@>\bigaw@\global\bigaw@\wdz@\fi
 \ifdim\wd@ne>\bigaw@\global\bigaw@\wd@ne\fi
 \ifCD@\enskip\fi
 \mathrel{\mathop{\hbox to\bigaw@{$%
  \setboxz@h{$\displaystyle-\m@th$}\ht\z@\z@
  \displaystyle\m@th\mathord\ifx\@twohead1\twoheadleftarrow\else\leftarrow\fi
  \mkern-6mu\cleaders
  \hbox{$\displaystyle\mkern-2mu\copy\z@\mkern-2mu$}\hfill
  \mkern-6mu\box\z@\ifx\@hook1\mkern-9mu\rhook\fi$}}%
 \ifdim\wd\tw@>\z@\limits^{#1}_{#2}\else\limits^{#1}\fi}%
 \ifCD@\enskip\fi\ampersand@}


\define\today{\number\day\ \ifcase\month\or
  January\or February\or March\or April\or May\or June\or
  July\or August\or September\or October\or November\or December\fi
  \ \number\year}
\def\pr@m@s{\ifx'\next\let\nxt\pr@@@s \else\ifx^\next\let\nxt\pr@@@t
  \else\let\nxt\egroup\fi\fi \nxt}

\define\widebar#1{\mathchoice
  {\setbox0\hbox{\mathsurround\z@$\displaystyle{#1}$}\dimen@.1\wd\z@
    \ifdim\wd\z@<.4em\relax \dimen@ -.16em\advance\dimen@.5\wd\z@ \fi
    \ifdim\wd\z@>2.5em\relax \dimen@.25em\relax \fi
    \kern\dimen@ \overline{\kern-\dimen@ \box0\kern-\dimen@}\kern\dimen@}%
  {\setbox0\hbox{\mathsurround\z@$\textstyle{#1}$}\dimen@.1\wd\z@
    \ifdim\wd\z@<.4em\relax \dimen@ -.16em\advance\dimen@.5\wd\z@ \fi
    \ifdim\wd\z@>2.5em\relax \dimen@.25em\relax \fi
    \kern\dimen@ \overline{\kern-\dimen@ \box0\kern-\dimen@}\kern\dimen@}%
  {\setbox0\hbox{\mathsurround\z@$\scriptstyle{#1}$}\dimen@.1\wd\z@
    \ifdim\wd\z@<.28em\relax \dimen@ -.112em\advance\dimen@.5\wd\z@ \fi
    \ifdim\wd\z@>1.75em\relax \dimen@.175em\relax \fi
    \kern\dimen@ \overline{\kern-\dimen@ \box0\kern-\dimen@}\kern\dimen@}%
  {\setbox0\hbox{\mathsurround\z@$\scriptscriptstyle{#1}$}\dimen@.1\wd\z@
    \ifdim\wd\z@<.2em\relax \dimen@ -.08em\advance\dimen@.5\wd\z@ \fi
    \ifdim\wd\z@>1.25em\relax \dimen@.125em\relax \fi
    \kern\dimen@ \overline{\kern-\dimen@ \box0\kern-\dimen@}\kern\dimen@}%
  }

\catcode`\@\active

\let\PVstyle=d 

\font\tenscr=rsfs10 
\font\sevenscr=rsfs7 
\font\fivescr=rsfs5 
\skewchar\tenscr='177 \skewchar\sevenscr='177 \skewchar\fivescr='177
\newfam\scrfam \textfont\scrfam=\tenscr \scriptfont\scrfam=\sevenscr
\scriptscriptfont\scrfam=\fivescr
\define\scr#1{{\fam\scrfam#1}}
\let\Cal\scr

\let\0\relax 
\define\boldProj{\operatorname{\text{\bc Proj}}}
\define\dist{\operatorname{dist}}
\define\exc{{\text{exc}}}
\define\Fl{\operatorname{Fl}}

\define\Gr{\operatorname{Gr}}
\define\ord{\operatorname{ord}}

\define\pr{\operatorname{pr}}
\define\Spec{\operatorname{Spec}}
\define\Supp{\operatorname{Supp}}
\define\restrictedto#1{\big|_{#1}}

\topmatter
\title On McQuillan's ``tautological inequality'' and the Weyl-Ahlfors
  theory of associated curves\endtitle
\rightheadtext{On McQuillan's inequality and the Weyl-Ahlfors theory}
\author Paul Vojta\endauthor
\affil University of California, Berkeley\endaffil
\address Department of Mathematics, University of California,
  970 Evans Hall\quad\#3840, Berkeley, CA \ 94720-3840\endaddress
\date \today \enddate
\thanks Supported by NSF grant DMS-0500512, and by MSRI.\endthanks
\subjclassyear{2000}
\subjclass Primary 32H30; Secondary 11J97, 14M15\endsubjclass

\abstract
In 1941, L. Ahlfors gave another proof of a 1933 theorem of H. Cartan on
approximation to hyperplanes of holomorphic curves in $\Bbb P^n$.
Ahlfors' proof built on earlier work of H. and J. Weyl (1938), and
proved Cartan's theorem by studying the associated curves of the holomorphic
curve.  This work has subsequently been reworked by H.-H. Wu in 1970,
using differential geometry, M. Cowen and P. A. Griffiths in 1976, further
emphasizing curvature, and by Y.-T. Siu in 1987 and 1990, emphasizing
meromorphic connections.  This paper gives another variation of the
proof, motivated by successive minima as in the proof of Schmidt's
Subspace Theorem, and using McQuillan's ``tautological inequality.''
In this proof, essentially all of the analysis is encapsulated within
a modified McQuillan-like inequality, so that most of the proof primarily
uses methods of algebraic geometry, in particular flag varieties.
A diophantine conjecture based on McQuillan's inequality is also posed.
\endabstract
\endtopmatter

\document

Cartan's theorem on value distribution with respect to hyperplanes
in $\Bbb P^n$ in general position \cite{Ca} remains to this day as
one of the pillars of value distribution of holomorphic curves in
higher-dimensional spaces.

\thm{\00.1} (Cartan)  Let $n\in\Bbb Z_{>0}$, and let $H_1,\dots,H_q$
be hyperplanes in $\Bbb P^n_{\Bbb C}$ in general position (i.e., all
collections of up to $n+1$ of the corresponding linear forms on $\Bbb C^{n+1}$
are linearly independent).  Let $f\:\Bbb C\to\Bbb P^n$ be a holomorphic
curve whose image is not contained in any hyperplane.  Then
$$\sum_{i=1}^q m_f(H_i,r)
  \le_\exc (n+1)T_f(r) + O(\log^{+}T_f(r)) + o(\log r)\;.$$
Here the notation $\le_\exc$ means that the inequality holds
for all $r$ outside a union of intervals of finite total length.
\endit

Cartan's proof is short and clever, but (so far) has not given much insight
into the analogies with number theory \cite{O}, \cite{V~1}.

A few years after Cartan's proof came out, H. and J. Weyl \cite{W-W}
and L. Ahlfors \cite{A} developed a different proof of this theorem,
based on a theory of associated curves modeled after the algebraic case.
Briefly, if $f\:\Bbb C\to\Bbb P^n(\Bbb C)$ is a holomorphic curve whose
image is not contained in any hyperplane, and if
$\bold x\:\Bbb C\to\Bbb C^{n+1}$ is a lifting of $f$, then for $d=1,\dots,n$
the $d^{\text{th}}$ associated curve of $f$ is the holomorphic curve
from $\Bbb C$ to the Grassmannian $\Gr_d(\Bbb C^{n+1})$ given for general
$z\in\Bbb C$ by the linear subspace spanned by the vectors
$\bold x(z),\bold x'(z),\dots,\bold x^{(d-1)}(z)$.  In contrast to Cartan's
proof, this proof is quite technical, and in fact Ahlfors notes in his
{\it Collected Works\/} that a reviewer had described it as a ``tour de force.''

This latter proof has been revisited over the years.  In 1970 H.-H. Wu
\cite{Wu} revisited the theory, revising it with an emphasis on differential
geometry.  In 1976 M. Cowen and
P. Griffiths \cite{Co-G} reworked the proof with further emphasis on curvature
(as part of a program of Griffiths to rework much of Nevanlinna theory
using curvature).  In the late 1980s, Y.-T. Siu reworked the proof again
using meromorphic connections \cite{S~1}, \cite{S~2}.

Although it is longer and more technical, Ahlfors' approach is better than
Cartan's from the point of view of comparison with the corresponding theorem
in number theory, Schmidt's Subspace Theorem.  Schmidt's theorem was
proved using the theory of successive minima, and many of the constructions
closely parallel those of associated curves.  Because of this, it is my
conviction that the derivative of a holomorphic function should translate
into number theory as some sort of object involving successive minima.
Most likely, the derivative of a holomorphic curve in a complex variety $X$,
modulo scalar multiplication (hence, a holomorphic curve in
$\Bbb P(\Omega_{X/\Bbb C})$, using Grothendieck's convention on
$\Bbb P(\Cal F)$ as described in Section \01) should correspond to the following
object in Arakelov theory.  Let $Y=\Spec\Cal O_k$ for a number field $k$
(or, let $Y$ be a smooth projective curve over a field of characteristic zero,
and let $k$ be its function field), let $X$ be an arithmetic variety over $Y$,
let $P\in X_k(k)$ be a $k$\snug-rational point, and let $i\:Y\to X$ be
the corresponding arithmetic curve in $X$.  Then the derivative should
correspond to a line subsheaf of maximal Arakelov degree in
$i^{*}\Omega_{X/\Cal O_k}\spcheck$.  See Section \02.

This comparison (except for some additional details) was explored in
Chapter 6 of \cite{V~1}.  That chapter presented Ahlfors' and Schmidt's
proofs of their respective theorems from the point of view of describing
their commonalities.  This effort was hamstrung, however, by a key difference
in their proofs.  Let $V$ be the vector space $\Bbb C^{n+1}$ or $k^{n+1}$
in the complex analytic and number theoretic contexts, respectively.
When Schmidt worked in $\bigwedge^d V$, he used hyperplanes corresponding
to $d$ of the original linear forms on $V$, hence an element of
$\bigwedge^d V\spcheck$.  Ahlfors, on the other hand, always worked with
a single hyperplane on $V$, via a more general interior product,
defined (in one special case) as follows.  Let $X\in\bigwedge^d V$
and let $\bold b\in V\spcheck$.  Then $(X\cdot\bold b)$ is the element of
$\bigwedge^{d-1}V$ characterized by the condition that
$$(X\cdot\bold b)\cdot Z = (X\cdot(\bold b\wedge Z))$$
for all $Z\in\bigwedge^{d-1}V\spcheck$.

More recently, M. McQuillan formulated a {\it tautological inequality\/}
(Theorem \01.2 below) which provides an elegant geometrical statement closely
related to geometric generalizations of the lemma on the logarithmic derivative
formulated by J. Noguchi \cite{N~1} and \cite{N~2}, S. Lu \cite{Lu},
P.-M. Wong \cite{Wo-S}, K. Yamanoi \cite{Y}, and others.  This significantly
cleared up the picture for the situation of Cartan's theorem (and
Ahlfors' proof thereof), and enabled me to translate Schmidt's proof
into the complex analytic case, obtaining a proof similar to Ahlfors',
but different in the sense that the associated curves in $\bigwedge^d V$
are now being compared with hyperplanes in $\bigwedge^d V\spcheck$;
i.e., the generalized interior product is no longer present.  (Therefore,
this proof could be regarded as being closer to Cartan's than Ahlfors.')

In addition, the proof in this paper is phrased almost entirely in terms
of algebraic geometry:  almost all of the analysis is encapsulated in
McQuillan's result.

Following \cite{V~3}, we actually prove a slightly stronger theorem than
Theorem \00.1:

\thm{\00.2}  Let $n\in\Bbb Z_{>0}$ and let $H_1,\dots,H_q$ be hyperplanes
in $\Bbb P^n$ (not necessarily in general position).  For $j=1,\dots,q$
let $\lambda_{H_j}$ be a corresponding Weil function; e.g.,
$$\lambda_{H_j}(\bold z)
  = -\frac12\log\frac{|L_j(z_0,\dots,z_n)|^2}{|z_0|^2+\dots+|z_n|^2}\;,$$
where $L_j$ is a linear form associated to $H_j$ and $[z_0:\dots:z_n]$
are homogeneous coordinates for $\bold z\in\Bbb P^n(\Bbb C)\setminus H_j$.
Let $f\:\Bbb C\to\Bbb P^n$ be a holomorphic curve whose image is not
contained in any hyperplane.  Then
$$\int_0^{2\pi}\max_J\sum_{i\in J}\lambda_{H_i}(f(re^{\sqrt{-1}\theta}))
    \frac{d\theta}{2\pi}
  \le_\exc (n+1)T_f(r) - N_W(r) + O(\log^{+}T_f(r)) + o(\log r)\;,\tag\00.2.1$$
where the maximum is taken over all subsets $J$ of $\{1,\dots,q\}$
for which the hyperplanes $H_i$, $i\in J$, lie in general position,
and $N_W(r)$ is the counting function for the zeroes of the Wronskian of
a choice of homogeneous coordinates of $f$ (with no common zeroes).
\endit

The paper is organized as follows.  Section \01 gives McQuillan's result.
A little additional work is needed because a slightly more general theorem
is needed here.  Section \02 introduces a diophantine conjecture motivated
by McQuillan's result and the strengthening proved in Section \01.
Section \03 contains a geometric discussion of
Grassmannian and flag varieties.  Section \04 discusses successive minima
in the context of Schmidt's theorem and proof, and uses this discussion
to motivate the general outline of the proof that follows.  Section \05
proves the main step of the proof in the case $d=1$, and Section \06
gives the general case.  The proof of Theorem \00.2 concludes in Section \07.

It would be more natural to phrase the proof in terms of jet spaces
over $\Bbb P^n$ instead of spaces $\bigwedge^d V\spcheck$, but attempts
to do so did not succeed.

Most of the standard notations of Nevanlinna theory (proximity function
$m_f(D,r)$\snug, counting function $N_f(D,r)$, and height (characteristic)
function $T_{f,\Cal L}(r)$) are the standard ones, as in for example
\cite{N~2, \S1} (plus $T_{f,\Cal L}(r)=T_{f,c_1(\Cal L)}(r)$).

\beginsection{\01}{McQuillan's ``Tautological Inequality''}

In this section let $X$ be a smooth compact complex variety and let $D$
be a normal crossings divisor (i.e., a divisor whose only singularities
are locally of the form $z_1\dots z_r=0$ for a suitable local holomorphic
coordinate system $z_1,\dots,z_n$ on $X$).  Throughout this paper, normal
crossings divisors are assumed to be effective and reduced.

If $\Cal E$ is a vector sheaf on $X$, then $\Bbb P(\Cal E)$ is defined as
$\boldProj\bigoplus_{d\ge0}S^d\Cal E$, so that points on $\Bbb P(\Cal E)$
lying over a point $x\in X$ correspond naturally to {\it hyperplanes\/} in
the fiber $\Cal E_x/\frak m_x\Cal E_x$ of $\Cal E$ at $x$.

Let $f\:\Bbb C\to X$ be a non-constant holomorphic curve whose image
is not contained in the support of $D$.  Then $f$ lifts to a holomorphic map
$f'\:\Bbb C\to\Bbb P(\Omega_{X/\Bbb C}(\log D))$.

\defn{\01.1}  The {\bc $D$\snug-modified ramification counting function}
of $f$ is the counting function for vanishing of the pull-back
$f^{*}\Omega_{X/\Bbb C}(\log D)$ (i.e., it counts the smallest order of
vanishing of $f^{*}s$ at $z\in\Bbb C$ for local sections $s$ of
$\Omega_{X/\Bbb C}(\log D)$ near $f(z)$).
\endit

McQuillan's tautological inequality is then the following.  It appeared
in \cite{McQ~1} with $D=0$ and \cite{McQ~2, V.1.2} in general.  See also
\cite{V~4, Thm.~A.6}.

\thm{\01.2}  Let $X$, $D$, $f$, and $f'$ be as above, and let $\Cal A$ be
a line sheaf on $X$ whose restriction to the Zariski closure of the image
of $f$ is big.  Then
$$T_{\Cal O(1),f'}(r) \le_\exc N^{(1)}_f(D,r) - N_{\text{Ram}(D),f}(r)
  + O(\log^{+}T_{\Cal A,f}(r)) + o(\log r)\;.$$
\endit

Let us compare Theorem \01.2 with the corresponding statement with $D=0$.
Of course $\Cal O(1)$ is different, since $\Bbb P(\Omega_{X/\Bbb C})$
is not isomorphic to $\Bbb P(\Omega_{X/\Bbb C}(\log D))$.  Also the term
$N^{(1)}_f(D,r)$ goes away when $D=0$.  The difference in heights
can be written in a canonical way as the height relative to a certain Cartier
divisor on the closure of the graph of the canonical birational map
$\Bbb P(\Omega_{X/\Bbb C}(\log D))\dashrightarrow\Bbb P(\Omega_{X/\Bbb C})$,
and the counting function associated to that divisor almost exactly cancels
the term $N^{(1)}_f(D,r)$ (there are also differences in the ramification
term).  This makes sense, since McQuillan's inequality is derived from
a version of the lemma on the logarithmic derivative, which concerns only
proximity functions.

This paper will need a variation of McQuillan's inequality which refrains
from moving things over into counting functions, and which considers a finite
list of normal crossings divisors $D_1,\dots,D_\ell$ in place of $D$.
These divisors may be chosen independently, so that the support of their
sum need not have normal crossings.  In order to state this modified
McQuillan inequality, we define a function on $\Bbb P(\Omega_{X/\Bbb C})$
as follows.

\lemma{\01.3}  Let $X$ and $D$ be as above.  Let $\Gamma$
be the closure of the graph of the canonical birational map
$\Bbb P(\Omega_{X/\Bbb C}(\log D))\dashrightarrow\Bbb P(\Omega_{X/\Bbb C})$,
let $p\:\Gamma\to\Bbb P(\Omega_{X/\Bbb C}(\log D))$ and
$q\:\Gamma\to\Bbb P(\Omega_{X/\Bbb C})$ be the projection morphisms,
and let $\Cal O_{\log}(1)$ and $\Cal O(1)$ denote the tautological
line sheaves on $\Bbb P(\Omega_{X/\Bbb C}(\log D))$ and
$\Bbb P(\Omega_{X/\Bbb C})$, respectively.  Then there is a
Cartier divisor $E$ on $\Gamma$ and a canonical isomorphism
$p^{*}\Cal O_{\log}(1)\cong q^{*}\Cal O(1)\otimes\Cal O(E)$.

Moreover, let $\lambda$ be a Weil function for $E$ on $\Gamma$, and
let $z_1,\dots,z_n$ be local coordinates on $X$ such that $D$ is locally
given by $z_1\dotsm z_r=0$.  Then $\lambda$ corresponds to a function
$\mu$ on $\Bbb P(\Omega_{X/\Bbb C})$ satisfying
$$\mu = -\frac12\log\frac{|dz_1|^2+\dots+|dz_n|^2}
  {|dz_1/z_1|^2+\dots+|dz_r/z_r|^2+|dz_{r+1}|^2+\dots+|dz_n|^2} + O(1)
  \tag\01.3.1$$
on compact subsets of the coordinate patch, minus the support of $D$.
\endit

\demo{Proof}  By canonicity, we may work locally on $X$, so let $z_1,\dots,z_n$
be local coordinates as above.  If $U=X\setminus\Supp D$, then
$\Omega_{X/\Bbb C}$ and $\Omega_{X/\Bbb C}(\log D)$ are canonically isomorphic
over $U$, so the above birational map gives an isomorphism
$\phi\:\pi_{\log}^{-1}(U)\overset\sim\to\to\pi^{-1}(U)$, where
$\pi_{\log}\:\Bbb P(\Omega_{X/\Bbb C}(\log D))\to X$ and
$\pi\:\Bbb P(\Omega_{X/\Bbb C})\to X$ are the natural projections.
This isomorphism naturally gives an isomorphism between
$\Cal O_{\log}(1)\restrictedto{\pi_{\log}^{-1}(U)}$ and
$\phi^{*}\Cal O(1)\restrictedto{\pi^{-1}(U)}$.  This shows that $E$ exists.

To construct the Weil function, let $z_1,\dots,z_n$ be a local coordinate
system as above, and let $s=a_1dz_1+\dots+a_ndz_n$ be a local section
of $\Omega_{X/\Bbb C}$, where $a_i$ are local functions on $X$.
Then, regarding $s$ as a section of $\Cal O(1)$ on $\Bbb P(\Omega_{X/\Bbb C})$,
its divisor has a Weil function of the form
$$-\frac12\log\frac{|a_1dz_1+\dots+a_ndz_n|^2}{|dz_1|^2+\dots+|dz_n|^2} + O(1)
  \;.$$
Doing the same for $s$ as a section of $\Cal O_{\log}(1)$ gives a Weil function
of the form
$$-\frac12\log\frac{|a_1dz_1+\dots+a_ndz_n|^2}
  {|dz_1/z_1|^2+\dots+|dz_r/z_r|^2+|dz_{r+1}|^2+\dots+|dz_n|^2} + O(1)\;.$$
Pulling these back to $\Gamma$ and subtracting then gives (\01.3.1).\qed
\enddemo

\remk{\01.4}  Since $\mu$ is bounded from below, the pull-back of $E$ to
the normalization of $\Gamma$ is effective.
\endit

\defn{\01.5} \cite{V~2, Def.~7.1} A {\bc generalized Weil function} on a
variety is a function on a dense open subset of the variety that pulls back
to a Weil function on some blowing-up of the variety.
\endit

The function $\mu$ in Lemma \01.3 is an example of a generalized Weil function
that is not a Weil function.  This function may be thought of as an approximate
archimedean version of the truncated counting function at $D$.  For example,
near smooth points of $D$, where we may assume that $D$ is locally given
by $z_1=0$, the graph $\Gamma$ in Lemma \01.3 is the blowing-up of
$\Bbb P(\Omega_{X/\Bbb C})$ at $z_1=dz_1=0$, and $\mu$ is the proximity
function of the strict transform of $D$ in that blowing-up.  One may be
tempted to denote it $m_f^{(1)}(D,r)$.

The version of McQuillan's inequality to be used here is the following.

\thm{\01.6}  Let $X$ be a smooth complete complex variety, and let
$D_1,\dots,D_\ell$ be normal crossings divisors on $X$ (whose sum need
not have normal crossings support).  Let $f\:\Bbb C\to X$ be a non-constant
holomorphic curve whose image is not contained in the support of $\sum D_i$,
let $f'\:\Bbb C\to\Bbb P(\Omega_{X/\Bbb C})$ be the lifting of $f$, and
let $\Cal A$ be a line sheaf on $X$ whose restriction to the Zariski
closure of the image of $f$ is big.  Let $\mu_1,\dots,\mu_\ell$ be generalized
Weil functions on $\Bbb P(\Omega_{X/\Bbb C})$ obtained from $D_1,\dots,D_\ell$,
respectively, as in Lemma \01.3.  Finally, let $N_{\text{Ram},f}(r)$ be the
counting function for the ramification of $f$.  Then
$$T_{\Cal O(1),f'}(r)
  + \int_0^{2\pi}\max_{1\le i\le\ell}
    \mu_i(f'(re^{\sqrt{-1}\theta}))\frac{d\theta}{2\pi}
  + N_{\text{Ram},f}(r)
  \le_\exc O(\log^{+}T_{\Cal A,f}(r)) + o(\log r)\;.\tag\01.6.1$$
\endit

\demo{Proof}  This proof is just a straightforward adaptation of the proof
from \cite{V~4, Appendix}, making obvious changes to accommodate the multiple
divisors.  We start with an enhanced version of the geometric lemma on
the logarithmic derivative due to P.-M. Wong \cite{Wo-S, Thm.~A3}:

\lemma{\01.6.2}  Let $X$, $D_1,\dots,D_\ell$, $f$, and $\Cal A$ be as in
Theorem \01.6.  Let $m\in\Bbb Z_{>0}$, and for each $i=1,\dots,\ell$
let $\|\cdot\|_i$ be a continuous pseudo jet metric on the jet space
$J^mX(-\log D_i)$.  Also let $j^m_{D_i}f\:\Bbb C\to J^mX(-\log D_i)$
denote the $m^{\text{th}}$ jet lifting of $f$.  Then
$$\int_0^{2\pi}\log^{+}\max_{1\le i\le\ell}
    \|j^m_{D_i}f(re^{\sqrt{-1}\theta})\|_i\frac{d\theta}{2\pi}
  \le_\exc O(\log^{+}T_{\Cal A,f}(r)) + o(\log r)\;.\tag\01.6.2.1$$
\endit

\demo{Proof}  This is proved by reducing to the special case $\ell=1$,
which is proved already \cite{Wo-S, Thm~A3}.  Let $\pi\:X'\to X$
be an embedded resolution of the set $\bigcup\Supp D_i$, let $D'$
be the normal crossings divisor on $X'$ lying over this set, let $\|\cdot\|'$
be a continuous jet metric on $J^mX'(-\log D')$, and let $g\:\Bbb C\to X'$
be the holomorphic lifting of $f$.  By compactness, we have
$$\|j^m_{D_i}f(z)\|_i \ll \|j^m_{D'}g(z)\|'$$
for all $z\in\Bbb C$ and all $i=1,\dots,\ell$, with a constant independent
of $z$.  Therefore
$$\split \int_0^{2\pi}\log^{+}\max_{1\le i\le\ell}
    \|j^m_{D_i}f(re^{\sqrt{-1}\theta})\|_i\frac{d\theta}{2\pi}
  &\le \int_0^{2\pi}\log^{+}
    \|j^m_{D'}g(re^{\sqrt{-1}\theta})\|'\frac{d\theta}{2\pi} + O(1) \\
  &\le_\exc O(\log^{+}T_{\pi^{*}\Cal A,g}(r)) + o(\log r) \\
  &\le O(\log^{+}T_{\Cal A,f}(r)) + o(\log r)\;, \endsplit$$
as was to be shown.\qed
\enddemo

We now recall some notation from \cite{V~4}.  Let $\Gamma$ be the closure
of the rational map $\Bbb P(\Omega_{X/\Bbb C}\oplus\Cal O_X)
 \dashrightarrow\Bbb P(\Omega_{X/\Bbb C})$, and let
$p\:\Gamma\to\Bbb P(\Omega_{X/\Bbb C}\oplus\Cal O_X)$
and $q\:\Gamma\to\Bbb P(\Omega_{X/\Bbb C})$ be the canonical projections.
The graph $\Gamma$ is obtained from $\Bbb P(\Omega_{X/\Bbb C}\oplus\Cal O_X)$
by blowing up the zero section of
$\Bbb V(\Omega_{X/\Bbb C})\subseteq\Bbb P(\Omega_{X/\Bbb C}\oplus\Cal O_X)$;
let $[0]$ denote the exceptional divisor.  The pull-backs of the tautological
line sheaves on $\Bbb P(\Omega_{X/\Bbb C}\oplus\Cal O_X)$ and
$\Bbb P(\Omega_{X/\Bbb C})$ are related by
$$q^{*}\Cal O(1) \cong p^{*}\Cal O(1)\otimes\Cal O(-[0])\;.\tag\01.6.3$$
Also let $[\infty]$ denote the divisor
$\Bbb P(\Omega_{X/\Bbb C}\oplus\Cal O_X)\setminus\Bbb V(\Omega_{X/\Bbb C})$.
Since $\Cal O([\infty])\cong\Cal O(1)$ on
$\Bbb P(\Omega_{X/\Bbb C}\oplus\Cal O_X)$ and since the lifted curve
$\partial f\:\Bbb C\to\Bbb P(\Omega_{X/\Bbb C}\oplus\Cal O_X)$ never meets
$[\infty]$, we have
$$\split T_{\Cal O(1),\partial f}(r) &= m_{\partial f}([\infty],r) + O(1) \\
  &= \int_0^{2\pi}\log^{+}\|Tf(re^{\sqrt{-1}\theta})\|\frac{d\theta}{2\pi}
    + O(1)\endsplit\tag\01.6.4$$
for a continuous metric $\|\cdot\|$ on the tangent bundle $TX$, which
we now fix.  Let $g\:\Bbb C\to\Gamma$ be the lifting of $f$ (so that
$p\circ g=\partial f$ and $q\circ g=f'$).  By definition,
$N_{\text{Ram},f}(r)=N_g([0],r)$, so (\01.6.3) and (\01.6.4) combine to give
$$T_{\Cal O(1),f'}(r) + N_{\text{Ram},f}(r)
  = \int_0^{2\pi}\log^{+}\|Tf(re^{\sqrt{-1}\theta})\|\frac{d\theta}{2\pi}
    - m_g([0],r) + O(1)\;.\tag\01.6.5$$

\lemma{\01.6.6}  Let $D$ be a normal crossings divisor on $X$; let $\mu$
be the corresponding generalized Weil function, as in Lemma \01.3; fix
continuous metrics $\|\cdot\|$ and $\|\cdot\|_D$ on $TX$ and $TX(-\log D)$,
respectively; and let $\lambda_{[0]}$ be a Weil function for the divisor
$[0]$ on $\Gamma$.  Then
$$\log^{+}\|Tf(z)\| + \mu(f'(z)) - \lambda_{[0]}(g(z))
  \le \log^{+}\|T_Df(z)\|_D + O(1)$$
for all $z\in\Bbb C$, where the constant in $O(1)$ is independent of $f$
and $z$.
\endit

\demo{Proof}  Choose local coordinates $z_1,\dots,z_n$ near $f(z)$ such that
$D$ is given locally by $z_1\dots z_r=0$.  Write $f=(f_1,\dots,f_n)$ in these
coordinates.  We then have
$$\align \log^{+}\|Tf(z)\|
    &= \frac12\log^{+}\bigl(|f_1'(z)|^2+\dots+|f_n'(z)|^2\bigr) + O(1)\;, \\
  \log^{+}\|T_Df(z)\|_D
    &= \frac12\log^{+}\biggl(\fracwithdelims||{f_1'(z)}{f_1(z)}^2
    + \dots + \fracwithdelims||{f_r'(z)}{f_r(z)}^2 \\
    &\qquad+ |f_{r+1}'(z)|^2 + \dots + |f_n'(z)|^2\biggr) + O(1)\;, \\
  \lambda_{[0]}(g(z)) &= \frac12\log^{+}\frac1{|f_1'(z)|^2+\dots+|f_n'(z)|^2}
    + O(1)\;, \\
  \intertext{and, from (\01.3.1),}
  \mu(f'(z)) &= -\frac12\log\bigl(|f_1'(z)|^2+\dots+|f_n'(z)|^2\bigr) \\
    &\qquad+ \frac12\log\biggl(\fracwithdelims||{f_1'(z)}{f_1(z)}^2
    + \dots + \fracwithdelims||{f_r'(z)}{f_r(z)}^2 \\
    &\qquad\qquad+ |f_{r+1}'(z)|^2 + \dots + |f_n'(z)|^2\biggr) + O(1)\;.
  \endalign$$
Therefore
$$\split \log^{+}\|T_Df(z)\|_D
  &\ge \frac12\log\biggl(\fracwithdelims||{f_1'(z)}{f_1(z)}^2
    + \dots + \fracwithdelims||{f_r'(z)}{f_r(z)}^2 \\
    &\qquad+ |f_{r+1}'(z)|^2 + \dots + |f_n'(z)|^2\biggr) + O(1) \\
  &= \mu(f'(z)) + \frac12\log\bigl(|f_1'(z)|^2+\dots+|f_n'(z)|^2\bigr) + O(1) \\
  &= \mu(f'(z)) + \log^{+}\|Tf(z)\| - \lambda_{[0]}(g(z)) + O(1)\;,
  \endsplit$$
as was to be shown.\qed
\enddemo

Applying this lemma with $D=D_i$ for $i=1,\dots,\ell$, taking the max,
and integrating gives
$$\split& \int_0^{2\pi}\log^{+}\|Tf(re^{\sqrt{-1}\theta})\|\frac{d\theta}{2\pi}
  + \int_0^{2\pi}\max_{1\le i\le\ell}
    \mu_i(f'(re^{\sqrt{-1}\theta}))\frac{d\theta}{2\pi}
    - m_g([0],r) \\
  &\qquad \le \int_0^{2\pi}\log^{+}\max_{1\le i\le\ell}
    \|T_{D_i}f(re^{\sqrt{-1}\theta})\|_i\frac{d\theta}{2\pi}\endsplit$$
and therefore
$$\split& T_{\Cal O(1),f'}(r) + N_{\text{Ram},f}(r)
  + \int_0^{2\pi}\max_{1\le i\le\ell}
    \mu_i(f'(re^{\sqrt{-1}\theta}))\frac{d\theta}{2\pi} \\
  &\qquad\le \int_0^{2\pi}\log^{+}\max_{1\le i\le\ell}
    \|T_{D_i}f(re^{\sqrt{-1}\theta})\|_i\frac{d\theta}{2\pi}\endsplit$$
by (\01.6.5).  Combining this with Lemma \01.6.2 (with $m=1$) then gives
(\01.6.1).\qed
\enddemo

\remk{\01.7}  One can also formulate conjectural Second Main Theorems with
multiple divisors $D_1,\dots,D_\ell$, and show that they would follow
from the $\ell=1$ case by the same methods as were used in the proof of
Lemma \01.6.2.
\endit

\remk{\01.8}  Theorem \01.6 can also be proved in the situation of a
finite ramified covering $p\:Y\to\Bbb C$ and a holomorphic curve $f\:Y\to X$.
In this case one would add a term $N_{\text{Ram},p}(r)$ on the right-hand
side of (\01.6.1).  Details are left to the reader.
\endit

\beginsection{\02}{A diophantine conjecture}

As noted in the Introduction, in the number field (or function field) case,
the derivative should be expressed in terms of successive minima, which
in Arakelov theory corresponds to to successive maxima (of degrees of
vector subsheaves).  We describe this in more detail as follows, giving
only the number field case since the translation to function fields is
straightforward.

Let $k$ be a number field, let $Y=\Spec\Cal O_k$, and let $\pi\:X\to Y$ be a
proper arithmetic variety.  Points $P\in X_k(k)$ correspond bijectively
to sections $i\:Y\to X$ of $\pi$.  As noted in the Introduction,
comparisons between Schmidt's proof of his Subspace Theorem and Ahlfors'
proof of Cartan's theorem suggest that derivatives of a holomorphic function
should correspond somehow to a line subsheaf in $i^{*}\Omega_{X/Y}\spcheck$
of maximal degree.  Of course $i^{*}\Omega_{X/Y}$ can only be
assumed to be a vector sheaf if $\pi$ is smooth, and this restriction is too
strong for most applications.  We can address this, however, as follows.
A line subsheaf of $i^{*}\Omega_{X/Y}\spcheck$ of maximal degree
corresponds to a quotient line sheaf of $i^{*}\Omega_{X/Y}$ of minimal degree,
and this corresponds to a lifting $i'\:Y\to\Bbb P(\Omega_{X/Y})$ of
$i\:Y\to X$ for which $(i')^{*}\Cal O(1)$ has minimal degree, where $\Cal O(1)$
is the tautological line sheaf.  Of course this degree is none other than
the height of the lifted rational point relative to $\Cal O(1)$.
This leads to the following conjecture, which should give a number-theoretic
analogue of Theorem \01.6.

\conj{\02.1}  Let $k$ and $Y$ be as above, let $S$ be a finite set of places
of $k$ containing all of the archimedean places, let $\pi\:X\to Y$ be a proper
arithmetic variety with smooth projective generic fiber $X_k$,
let $D_1,\dots,D_\ell$ be effective Cartier divisors on $X$ whose restrictions
to $X_k$ are normal crossings divisors, and let $\mu_1,\dots,\mu_\ell$ be
generalized Weil functions on $X_k$ obtained from $D_1,\dots,D_\ell$,
respectively, as in Lemma \01.3.  Fix an ample line sheaf $\Cal A$ on $X_k$,
fix absolute heights $h_{\Cal O(1)}$ on $\Bbb P(\Omega_{X/Y})$ and $h_{\Cal A}$
on $X$, and fix constants $\epsilon>0$ and $C$.  For $P\in X_k(\bar k)$
let $S_P$ be the set of places of $k(P)$ lying over places in $S$.
Then, for all but finitely many $P\in X_k(\bar k)$ not lying in the
support of any $D_i$, there exists a point $P'\in\Bbb P(\Omega_{X/Y})(k(P))$
lying over $P$ and satisfying the inequality
$$h_{\Cal O(1)}(P')
    + \frac1{[k(P):\Bbb Q]}\sum_{w\in S_P}\max_{1\le i\le\ell}\mu_{i,w}(P')
  \le d(P) + \epsilon\,h_{\Cal A}(P) + C\;.\tag\02.1.1$$
\endit

Here
$$d(P)=\frac1{[k(P):\Bbb Q]}\log|D_{k(P)}|$$
is as in \cite{V~1, pp.~57--58}.  Of course, if we restrict to rational
points $P\in X_k(k)$ then this term goes away.  Also, if $\ell=0$ then
(\02.1.1) becomes further shortened to
$$h_{\Cal O(1)}(P') \le \epsilon\,h_{\Cal A}(P) + O(1)\;.\tag\02.2$$

The restriction that $X_k$ be projective is not essential, but eliminating
it would require an adequate replacement for $h_{\Cal A}$, which would take
some work.

The point $P'$ may not be uniquely determined by (\02.1.1); however, Schmidt's
proof focuses on points where there is a gap in the successive minima, and
such a gap would cause $P'$ to be uniquely defined.  One could then drop
the requirement that $P'$ be rational over $k(P)$.

This conjecture is somewhat reminiscent of Szpiro's work on small points.
Some differences include the (non-essential) fact that the above conjecture
looks only at points rational over $k(P)$, the fact that this conjecture
is phrased in a relative setting, and the fact that we are looking at the
smallest height rather than a lim inf of heights.

This conjecture will be developed further in a subsequent paper.

\beginsection{\03}{Some Geometry of Grassmann and Flag Varieties}

This section provides a basic result about flag and Grassmann varieties
(Proposition \03.7).  Although we only need the result in the case
$X=\Spec\Bbb C$, we work over an arbitary scheme $X$ since it is not
any harder, and it may be useful in later work.

\narrowthing{}  Throughout this section, $\Cal E$ is a vector sheaf
of rank $r$ over a scheme $X$.
\endit

The basic idea is as follows.  Let $s\in\Bbb Z_{>0}$ and let
$r\ge d_1>\dots>d_s\ge 0$ be integers.  The flag variety is a scheme over
$X$ whose fiber over $x\in X$ parametrizes flags of linear subspaces
$$W_1\subseteq W_2\subseteq \dots \subseteq W_s\subseteq \Cal E(x)\;,$$
where $\Cal E(x)$ is the fiber of $\Cal E$ over $x$, and where $W_i$
is a $k(x)$\snug-vector subspace of the fiber $\Cal E(x)$ of codimension $d_i$
for all $i$.  However, it is useful to phrase this definition as a moduli
problem.

\defn{\03.1}  Let $s\in\Bbb Z_{>0}$ and let $r\ge d_1>\dots>d_s\ge 0$
be integers.  Then the {\bc flag bundle} $\Fl^{d_1,\dots,d_s}(\Cal E)$
(if it exists) is the $X$\snug-scheme $\pi\:\Fl^{d_1,\dots,d_s}(\Cal E)\to X$
representing the contravariant functor $F$ from $X$\snug-schemes to sets,
defined as follows.  If $\phi\:T\to X$ is an $X$\snug-scheme, then $F(T)$
is the set of all flags
$$\Cal F_1\subseteq\dots\subseteq\Cal F_s\tag\03.1.1$$
of vector subsheaves of $\phi^{*}\Cal E$ such that $\phi^{*}\Cal E/\Cal F_i$
is a vector sheaf of rank $d_i$ for all $i$.  In the special case $s=1$,
the flag bundle $\Fl^d(\Cal E)$ (for $0\le d\le r$) is also called the
{\bc Grassmann bundle}, and is denoted $\Gr^d(\Cal E)$.  It represents
the functor of vector subsheaves whose quotient is a vector sheaf of rank $d$.
\endit

\remk{\03.2}  Often it is more convenient to refer to chains of surjections
$$\phi^{*}\Cal E \twoheadrightarrow \Cal G_s
  \twoheadrightarrow\dots\twoheadrightarrow \Cal G_1\;,$$
where each $\Cal G_i$ is a vector sheaf of rank $d_i$, in place of flags
(\03.1.1).  The flag bundle also represents the functor of these objects,
up to an obvious notion of isomorphism.
\endit

As is customary with Chow groups, superscripts indicate codimension and
subscripts indicate dimension, so we also write
$$\Fl_{d_1,\dots,d_s}(\Cal E) = \Fl^{r-d_1,\dots,r-d_s}(\Cal E)
  \qquad\text{and}\qquad
  \Gr_d(\Cal E) = \Gr^{r-d}(\Cal E)\;,$$
where $0\le d_1<\dots<d_s\le r$ and $0\le d\le r$, respectively.

Finally, if $X=\Spec k$ for a field $k$, and $V$ is a finite-dimensional
vector space over $k$, then we write
$\Fl^{d_1,\dots,d_s}(V)=\Fl^{d_1,\dots,d_s}(\widetilde V)$ and
$\Gr^d(V)=\Gr^d(\widetilde V)$, etc.  These are the {\bc flag varieties}
and {\bc Grassmann varieties}, respectively.

\prop{\03.3}  Let $s\in\Bbb Z_{>0}$ and let $r\ge d_1>\dots>d_s\ge 0$
be integers.  Then the flag bundle $\pi\:\Fl^{d_1,\dots,d_s}(\Cal E)\to X$
exists, and has universal sheaves
$$\Cal U_1\subseteq\dots\subseteq\Cal U_s\subseteq\pi^{*}\Cal E$$
such that if $f\:T\to\Fl^{d_1,\dots,d_s}(\Cal E)$ corresponds to a flag
(\03.1.1), then $\Cal F_i=f^{*}\Cal U_i$ for all $i$.
\endit

\demo{Proof}  Following \cite{F~2, Prop.~14.2.1}, we use induction on $s$.
The case $s=1$ (the Grassmannian) was already proved by Kleiman
\cite{K, Prop.~1.2}.  If $s>1$ then let
$\Cal U_2'\subseteq\dots\subseteq\Cal U_s'$ be the universal sheaves
on $\Fl^{d_2,\dots,d_s}(\Cal E)$; then
$$\Fl^{d_1,\dots,d_s}(\Cal E)=\Gr^{d_2-d_1}(\Cal U_2')$$
represents the functor of Definition \03.1, with universal sheaf $\Cal U_1$
equal to the universal sheaf of the Grassmannian and $\Cal U_i$ equal
to the pull-backs of $\Cal U_i'$ for $i=2,\dots,s$.\qed
\enddemo

We note that $\Gr^1(\Cal E)=\Bbb P(\Cal E)$, and that its universal sheaf
is the kernel of the canonical map $\pi^{*}\Cal E\to\Cal O(1)$
\cite{EGA, II~4.2.5}.

If $d_1,\dots,d_s$ are as above, and if $t\in\Bbb Z_{>0}$ and
$r\ge e_1>\dots>e_t\ge0$ are integers such that $\{e_1,\dots,e_t\}$ is
a subset of $\{d_1,\dots,d_s\}$, then there is a {\bc forgetful morphism}
$$\Fl^{d_1,\dots,d_s}(\Cal E) @>>> \Fl^{e_1,\dots,e_t}(\Cal E)$$
over $X$, and those universal bundles $\Cal U_i$ on
$\Fl^{d_1,\dots,d_s}(\Cal E)$ for which $d_i=e_j$ for some $j$, are
pull-backs of the corresponding universal bundles $\Cal U_j'$ on
$\Fl^{e_1,\dots,e_t}(\Cal E)$.  If $0\le i<j\le s$, then the diagram
$$\CD \Fl^{d_1,\dots,d_s}(\Cal E) @>>> \Fl^{d_{i+1},\dots,d_s}(\Cal E) \\
  @VVV @VVV \\
  \Fl^{d_1,\dots,d_j}(\Cal E) @>>> \Fl^{d_{i+1},\dots,d_j}(\Cal E)\endCD
  \tag\03.4$$
in which the arrows are forgetful morphisms, is Cartesian.  Indeed, this is
clear from the definition of flag varieties as objects representing functors,
since if $\Cal F_{i+1},\dots,\Cal F_j$ in (\03.1.1) are fixed, then
$\Cal F_1,\dots,\Cal F_i$ can be chosen independently of
$\Cal F_{j+1},\dots,\Cal F_s$.

We now look at line sheaves on flag bundles.  For the case of flag varieties,
a reference is \cite{F~1, Ch.~9}.

For $0\le d\le r$ we define a morphism
$i\:\Gr^d(\Cal E)\to\Bbb P\bigl(\bigwedge^d\Cal E\bigr)$, by specifying
(via Remark \03.2) that a $T$\snug-point corresponding to a surjection
$\phi^{*}\Cal E\twoheadrightarrow\Cal G$ is taken to the natural surjection
$\phi^{*}\bigwedge^d\Cal E\twoheadrightarrow\bigwedge^d\Cal G$;
here $\phi\:T\to X$ is the structural morphism.  By \cite{K, Prop.~1.5}
this is a closed immersion.  It is called the {\bc Pl\"ucker embedding}.
From this definition (in particular with $T=\Gr^d(\Cal E)$) it follows that
$i^{*}\Cal O(1)=\bigwedge^d(\pi^{*}\Cal E/\Cal U)$, where $\Cal U$ is
the universal sheaf on $\Gr^d(\Cal E)$.  This line sheaf $i^{*}\Cal O(1)$
is also denoted $\Cal O(1)$.

The Grassmann coordinates induce a cover of
$\Gr^d\bigl(\Cal E\restrictedto U\bigr)$ by open subsets isomorphic
to $\Bbb A^{d(r-d)}_U$, for all open $U$ on $X$ for which
$\Cal E\restrictedto U$ is trivial \cite{K, Prop.~1.6}.
Therefore $\Gr^d(\Cal E)$ is smooth over $X$ with connected fibers.
By induction on $s$, the same holds for flag bundles.

Returning to flag bundles, let $r\ge d_1>\dots>d_s\ge0$, let
$$\pr_{d_i}\:\Fl^{d_1,\dots,d_s}(\Cal E)\to\Gr^{d_i}(\Cal E)$$
denote the forgetful morphisms, and for all $a_1,\dots,a_s\in\Bbb Z$
define the line sheaf
$$\split \Cal O(a_1,\dots,a_s)
  &= \pr_{d_1}^{*}\Cal O(a_1)\otimes\dots\otimes\pr_{d_s}^{*}\Cal O(a_s) \\
  &\cong
    \Bigl(\bigwedge\nolimits^{d_1}(\pi^{*}\Cal E/\Cal U_1)\Bigr)^{\otimes a_1}
    \otimes\dots\otimes
    \Bigl(\bigwedge\nolimits^{d_s}(\pi^{*}\Cal E/\Cal U_s)\Bigr)^{\otimes a_s}
  \endsplit\tag\03.5$$
on $\Fl^{d_1,\dots,d_s}(\Cal E)$, where $\Cal U_i$ are the universal subsheaves.

\lemma{\03.6}  Let $0<d<r$, let $G=\Gr^d(\Cal E)$, let $\pi\:G\to X$ be its
structural map, let $\Cal F\subseteq\pi^{*}\Cal E$ be the universal subbundle,
and let $\Cal Q=\pi^{*}\Cal E/\Cal F$ be the quotient.  Then
\roster
\myitem a.  There is a canonical isomorphism
$$\alpha\:\Fl^{d+1,d}(\Cal E) \overset\sim\to\to \Bbb P(\Cal F)$$
over $G$, and a canonical isomorphism
$$\alpha^{*}\Cal O(1) \cong \Cal O(1,-1)\;.\tag\03.6.1$$
\myitem b.  There is a canonical isomorphism
$$\beta\:\Fl^{d,d-1}(\Cal E) \overset\sim\to\to \Bbb P(\Cal Q\spcheck)$$
over $G$, and a canonical isomorphism
$$\beta^{*}\Cal O(1) \cong \Cal O(-1,1)\;.\tag\03.6.2$$
\endroster
\endit

\demo{Proof}  The existence of $\alpha$ follows from the fact that
$$\Bbb P(\Cal F) = \Gr^1(\Cal F) = \Fl^{d+1,d}(\Cal E)\;.$$
Let $\pi'\:\Fl^{d+1,d}(\Cal E)\to X$ be the canonical map, and let
$\Cal F_1\subseteq\Cal F_2$ be the universal bundles on $\Fl^{d+1,d}(\Cal E)$.
Applying \cite{H, II Ex.~5.16d} to the short exact sequence
$$0 @>>> \Cal F_2/\Cal F_1 @>>> (\pi')^{*}\Cal E/\Cal F_1
  @>>> (\pi')^{*}\Cal E/\Cal F_2 @>>> 0\tag\03.6.3$$
gives a canonical isomorphism
$$\bigwedge\nolimits^{d+1}\bigl((\pi')^{*}\Cal E/\Cal F_1\bigr)
  \cong (\Cal F_2/\Cal F_1)
    \otimes\bigwedge\nolimits^d\bigl((\pi')^{*}\Cal E/\Cal F_2\bigr)\;.$$
Now let $\rho\:\Bbb P(\Cal F)\to G$ be the structural map; then
$\alpha^{*}\rho^{*}\Cal F=\Cal F_2$ (as subbundles of $(\pi')^{*}\Cal E$), and
$\Cal F_1\cong\alpha^{*}\bigl(\ker\bigl(\rho^{*}\Cal F\to\Cal O(1)\bigr)\bigr)$,
so we have canonical isomorphisms
$$\split \alpha^{*}\Cal O(1) &\cong \Cal F_2/\Cal F_1 \\
  &\cong \bigwedge\nolimits^{d+1}\bigl((\pi')^{*}\Cal E/\Cal F_1\bigr)
    \otimes \bigwedge\nolimits^d\bigl((\pi')^{*}\Cal E/\Cal F_2\bigr)\spcheck \\
  &\cong \Cal O(1,-1)\endsplit$$
by (\03.5).  This gives (\03.6.1).

Dually, the existence of $\beta$ follows from the fact that
$$\Bbb P(\Cal Q\spcheck) = \Gr_1(\Cal Q) = \Fl^{d,d-1}(\Cal E)\;.$$
As before let $\pi'\:\Fl^{d,d-1}(\Cal E)\to X$ be the canonical map and
$\Cal F_1\subseteq\Cal F_2$ the universal subbundles.  The short exact
sequence (\03.6.3) now gives a canonical isomorphism
$$\bigwedge\nolimits^d\bigl((\pi')^{*}\Cal E/\Cal F_1\bigr)
  \cong (\Cal F_2/\Cal F_1)
    \otimes\bigwedge\nolimits^{d-1}\bigl((\pi')^{*}\Cal E/\Cal F_2\bigr)\;.$$
As in the first half of the proof, let $\rho\:\Bbb P(\Cal Q\spcheck)\to X$
be the structural map; then we have a natural isomorphism
$\beta^{*}\rho^{*}\Cal Q\spcheck
 \cong\bigl((\pi')^{*}\Cal E/\Cal F_1\bigr)\spcheck$, by which the subbundles
$\beta^{*}\ker\bigl(\rho^{*}\Cal Q\spcheck\to\Cal O(1)\bigr)$ and
$\bigl((\pi')^{*}\Cal E/\Cal F_2\bigr)\spcheck$ correspond.  Therefore,
$$\split \beta^{*}\Cal O(1)
  &\cong \bigl((\pi')^{*}\Cal E/\Cal F_1\bigr)\spcheck
    / \bigl((\pi')^{*}\Cal E/\Cal F_2\bigr)\spcheck \\
  &\cong (\Cal F_2/\Cal F_1)\spcheck \\
  &\cong \bigwedge\nolimits^{d-1}\bigl((\pi')^{*}\Cal E/\Cal F_2\bigr)
    \otimes \bigwedge\nolimits^d\bigl((\pi')^{*}\Cal E/\Cal F_1\bigr)\spcheck \\
  &\cong \Cal O(-1,1)\endsplit$$
by (\03.5) again, which gives (\03.6.2).\qed
\enddemo

\prop{\03.7}  Let $d$, $G$, $\pi$, $\Cal F$, and $\Cal Q$ be as in Lemma \03.6.
Then there is a canonical map
$$\phi\:\Fl^{d+1,d,d-1}(\Cal E) @>>> \Bbb P(\Omega_{G/X})$$
and a canonical isomorphism
$$\phi^{*}\Cal O(1) \cong \Cal O(1,-2,1)\;.$$
\endit

\demo{Proof}  By \cite{H, II Prop.~7.12}, it suffices to give a natural
surjection
$$\pr_d^{*}\Omega_{G/X}\to\Cal O(1,-2,1)\;,$$
where we recall that $\pr_d\:\Fl^{d+1,d,d-1}(\Cal E)\to G$ is the forgetful map.
By \cite{F~2, B.5.8}, $\Omega_{G/X}\cong\Cal F\otimes\Cal Q\spcheck$
(this was stated only for Grassmannians over fields, but the same proof is
valid for general $X$).  
By (\03.4) and Lemma \03.6,
$$\split \Fl^{d+1,d,d-1}(\Cal E)
  &\cong \Fl^{d,d-1}(\Cal E) \times_G \Fl^{d+1,d}(\Cal E) \\
  &\cong \Bbb P(\Cal F) \times_G \Bbb P(\Cal Q\spcheck)\;.\endsplit$$
Letting $p\:\Fl^{d+1,d,d-1}(\Cal E)\to\Bbb P(\Cal F)$ and
$q\:\Fl^{d+1,d,d-1}(\Cal E)\to\Bbb P(\Cal Q\spcheck)$ be the projection
morphisms via the above isomorphism, we have a canonical surjective map
$$\split \pr_d^{*}\Omega_{G/X}
  &\cong \pr_d^{*}\Cal F\otimes\pr_d^{*}\Cal Q\spcheck \\
  &@[twohead]>>> p^{*}\Cal O(1)\otimes q^{*}\Cal O(1) \\
  &\cong \Cal O(1,-1,0)\otimes\Cal O(0,-1,1) \\
  &\cong \Cal O(1,-2,1)\endsplit$$
by (\03.6.1) and (\03.6.2), as was to be shown.\qed
\enddemo

\cor{\03.8}  Let $V$ be a complex vector space of dimension $r$,
and let $P=\Bbb P(V\spcheck)=\allowmathbreak\Gr_1(V)$.  Then there is
a canonical isomorphism
$$\phi\:\Fl_{1,2}(V) @>\sim>> \Bbb P(\Omega_{P/\Bbb C})$$
of schemes over $P$ and a canonical isomorphism
$$\phi^{*}\Cal O(1) \cong \Cal O(-2,1)$$
of line sheaves on $\Fl_{1,2}(V)$.
\endit

\demo{Proof}  Apply Proposition \03.7 with $X=\Spec\Bbb C$,
$\Cal E=\widetilde V$, $d=r-1$, and note that there is a canonical isomorphism
$$\Fl^{r,r-1,r-2}(V) = \Fl_{0,1,2}(V) \cong \Fl_{1,2}(V)$$
taking $\Cal O(1,-2,1)$ to $\Cal O(-2,1)$ (note that $\Cal O(1,0,0)$ is
trivial on $\Fl_{0,1,2}(V)$).  This gives the map $\phi$ and an isomorphism
of line sheaves.

The fact that $\phi$ is an isomorphism is left to the reader, since it will
not be used in the sequel.\qed
\enddemo

\remk{\03.9}  The intuition behind Proposition \03.7 is as follows.
Let $X=\Spec\Bbb C$ and $\Cal E=\widetilde V$ for a finite-dimensional
complex vector space $V$.  Let $G=\Gr_d(V)$.  Giving a closed point in
$\Bbb P(\Omega_{G/\Bbb C})$ corresponds to giving a line in the Pl\"ucker
embedding of $G$ and a point on that line.  One can show that
this line is contained in the image of $G$ if and only if the deformation
corresponding to the tangent direction in $G$ at the point corresponds to
deforming the $d$\snug-dimensional subspace in such a way that some
$(d-1)$\snug-dimensional subspace within the original subspace also lies
within the deformations.  If so, then all these deformations lie within
some fixed $(d+1)$\snug-dimensional subspace.  This is exactly the situation
of the associated holomorphic curves $X^d$ appearing in Section \06.
\endit

\remk{\03.10}  Proposition \03.7 can be used to prove the Pl\"ucker formulas
\cite{G-H, p.~270}.  We also recall that \cite{W-W} described its results
as Pl\"ucker formulas for holomorphic curves.
\endit

\beginsection{\04}{Motivation from Successive Minima}

We start the main sequence of the proof of Theorem \00.2 by showing how
the main step is
motivated by the use of successive minima in Schmidt's proof.  First consider
the simplest case, over $\Bbb Q$ and considering only the infinite place.
Let $L_0,\dots,L_n$ be linearly independent linear forms in $n+1$ variables
$x_0,\dots,x_n$ with algebraic coefficients.  Let $\bold x$ be a point
close to the corresponding hyperplanes, and let $[x_0:\dots:x_n]$ be
homogeneous coordinates for $\bold x$.  We may assume that the $x_i$ are
integers and that they are (collectively) relatively prime.
Then the multiplicative and logarithmic heights of $\bold x$ are
$|\bold x|=\max_{0\le i\le n}|x_i|$ and $h(\bold x)=\log|\bold x|$,
respectively, and the proximity functions relative to the hyperplanes $L_i=0$
can be taken to be
$$m(L_i,\bold x) = -\log\frac{|L_i(\bold x)|}{\max_j|x_j|}
  = -\log|L_i(\bold x)| + h(\bold x)\;.$$

Since there are $n+1$ linear forms, Schmidt's inequality
$$\sum m(L_i,\bold x)\le (n+1+\epsilon)h(\bold x)+O(1)$$
reduces to
$\sum -\log|L_i(\bold x)| \le \epsilon\,h(\bold x) + O(1)$, or equivalently
$$\prod|L_i(\bold x)| \gg |\bold x|^\epsilon\;.\tag\04.1$$
Therefore it can be addressed by studying the successive minima of the
parallelepiped
$$A_i|L_i(\bold x)|\le 1,\qquad 0\le i\le n$$
in $\Bbb R^{n+1}$ for suitable positive real constants $A_0,\dots,A_n$.
Recall that the successive minima of this parallelepiped are defined by
the condition that the $d^{\text{th}}$ successive minimum is the smallest
real number $\lambda_d$ such that the scaled parallelepiped
$A_i|L_i(\bold x)|\le\lambda_d$ contains at least $d$ linearly independent
lattice points in $\Bbb Z^{n+1}$.

The proof proceeds by contradiction, so assume that there are infinitely
many counterexamples $\bold x$ to (\04.1).
The constants $A_i$ are chosen, depending on $\bold x$, such that
$A_0\dotsm A_n=1$, and so that the first successive minimum is small.
Choosing the $A_i$ so that $\bold x$ lands in the corner of the parallelepiped
gives the best upper bound for $\lambda_1$, so we require
$A_i|L_i(\bold x)| = \lambda_1$ for all $i$, resulting in
$$A_i = \frac{\root{n+1}\of{\prod_j|L_j(\bold x)|}}{|L_i(\bold x)|}
  \qquad\text{and}\qquad
  \lambda_1 = \root{n+1}\of{\prod_j|L_j(\bold x)|}\;.$$
Thus, assuming that (\04.1) fails for $\bold x$, this would imply
$\lambda_1\ll |\bold x|^{-\epsilon}$ (for a different $\epsilon$).
(There is a slight complication in case that some point other than $\bold x$
gives rise to the first successive minimum, but that need not concern us here.)

The theory of successive minima implies that $\lambda_1\dotsm\lambda_{n+1}$
is bounded away from $0$ by a constant depending only on $n$.  Schmidt's
proof derives a contradiction by showing that
$$\frac{\lambda_{d+1}}{\lambda_d} \ll |\bold x|^\epsilon$$
for all $1\le d\le n$ (and for a different $\epsilon>0$).

Translating the expression for $\lambda_1$ into Nevanlinna notation gives
$$-\log\lambda_1 = \frac{-\sum\log|L_i(\bold x)|}{n+1}
  = \frac{\sum\bigl(m(L_i,\bold x) - h(\bold x)\bigr)}{n+1}
  = \frac{\sum m(L_i,\bold x)}{n+1} - h(\bold x)\;.$$

For higher successive minima, say for $\lambda_d$ with $1\le d\le n+1$,
we use the fact that the product $\lambda_1\dotsm\lambda_d$ arises
(up to a bounded constant factor) as the first successive minimum in
$\bigwedge^d\Bbb R^{n+1}$ for the parallelepiped
$$A_{i_1}\dotsm A_{i_d}L_I(X)\le 1\;,$$
where $I=\{i_1,\dots,i_d\}$ varies over all $d$\snug-element subsets of
$\{0,\dots,n\}$ and
$$L_I=L_{i_1}\wedge\dots\wedge L_{i_d}\;.$$
If $\bold x,\bold x',\dots,\bold x^{(d-1)}$ are lattice points in $\Bbb Z^{n+1}$
corresponding to the first $d$ successive minima, and if
$X^d=\bold x\wedge\dots\wedge\bold x^{(d-1)}$, then $X^d$ gives rise to
the first successive minimum in $\bigwedge^d\Bbb R^{n+1}$ (or it comes within
a constant factor of doing so), so for $1\le d\le n+1$ (and also for $d=0$)
we have
$$\lambda_1\dotsm\lambda_d \gg\ll
  \max_{\#I=d}\frac{|L_I(X^d)|}{|L_{i_1}(\bold x)|\dotsm|L_{i_d}(\bold x)|}
    \lambda_1^d\;.$$
Therefore the ratio of consecutive successive minima can be obtained as
$$\split \frac{\lambda_{d+1}}{\lambda_d}
  &= \frac{(\lambda_1\dotsm\lambda_{d-1})(\lambda_1\dotsm\lambda_{d+1})}
    {(\lambda_1\dotsm\lambda_d)^2} \\
  &\gg\ll \frac{\max\limits_{\#I=d-1}
        \frac{|L_I(X^{d-1})|}{|L_{i_1}(\bold x)|\dotsm|L_{i_{d-1}}(\bold x)|}
      \cdot \max\limits_{\#I=d+1}
        \frac{|L_I(X^{d+1})|}{|L_{i_1}(\bold x)|\dotsm|L_{i_{d+1}}(\bold x)|}}
      {\Bigl(\max\limits_{\#I=d}
        \frac{|L_I(X^{d})|}{|L_{i_1}(\bold x)|\dotsm|L_{i_d}(\bold x)|}\Bigr)^2}
  \endsplit$$
for all $1\le d\le n$.

One would hope that the three maxima in the above expression occur at
subsets that are nested in a reasonable way, but there is no obvious reason
for why this should hold.  Therefore, we will deviate a little from Schmidt's
proof by taking geometric means instead of maxima.  So, in motivating the
structure of the proof of Theorem \00.2, we assume that
$$\split \lambda_1\dotsm\lambda_d
  &= \left(\prod_{\#I=d}
    \frac{|L_I(X^d)|}{|L_{i_1}(\bold x)|\dotsm|L_{i_d}(\bold x)|}\lambda_1^d
    \right)^{1\bigm/\binom{n+1}{d}} \\
  &= \lambda_1^d
      \frac{\left(\prod_{\#I=d}|L_I(X^d)|\right)^{1\bigm/\binom{n+1}{d}}}
      {\left(\prod_{i=0}^n |L_i(\bold x)|\right)^{\frac d{n+1}}} \\
  &= \left(\prod_{\#I=d}|L_I(X^d)|\right)^{1\bigm/\binom{n+1}{d}}\;.\endsplit$$
Then
$$\frac{\lambda_{d+1}}{\lambda_d}
  = \frac{\left(\prod_{\#I=d-1}|L_I(X^{d-1})|\right)^{1\bigm/\binom{n+1}{d-1}}
      \left(\prod_{\#I=d+1}|L_I(X^{d+1})|\right)^{1\bigm/\binom{n+1}{d+1}}}
    {\left(\prod_{\#I=d}|L_I(X^d)|\right)^{2\bigm/\binom{n+1}{d}}}\;.\tag\04.2$$
One would then bound this by $H(\bold x)^\epsilon$.

This derivation assumed that there was only one infinite place (and exactly
$n+1$ hyperplanes---the fact that Schmidt's theorem in this case is not
a trivial consequence of the First Main Theorem stems from the fact that
the coefficients are allowed to be algebraic).  In the general case, we work
with a collection of $n+1$ hyperplanes that is allowed to vary with the
place (or, in the Nevanlinna case, to vary with $z\in\Bbb C$ with finitely
many possibilities).

\beginsection{\05}{First Step of the Proof}

The main step in the proof is Proposition \06.2.  Here we prove a special case.
We start with some notation.

Let $V=\Bbb C^{n+1}$, with the usual norm $|\bold v|^2=|v_0|^2+\dots+|v_n|^2$
for $\bold v=(v_0,\dots,v_n)$ in $V$.  Let $\bold y\:\Bbb C\to V$ be a
holomorphic map, not identically zero, with coordinate functions
$y_0,\dots,y_n$.  Define
$$\widebar T_{\bold y}(r)
  = \int_0^{2\pi}\log\bigl|\bold y(re^{\sqrt{-1}\theta})\bigr|
    \frac{d\theta}{2\pi}\;,$$
and let $N_{\bold y}(r)$ be the counting function for simultaneous vanishing of
the coordinates $y_0,\dots,y_n$ of $\bold y$.

Write $P=\Bbb P(V\spcheck)$, so that $P\cong\Bbb P^n$ and $\bold y$
corresponds to a holomorphic curve $f\:\Bbb C\to\Bbb P^n$.  Note that
$$T_{\Cal O(1),f}(r) = \widebar T_{\bold y}(r) - N_{\bold y}(r)\;.$$
We also write $T_f(r)=T_{\Cal O(1),f}(r)$.

Assume also that $f$ is not constant.  Then $\bold y\wedge\bold y'$ is a
holomorphic map $\Bbb C\to\bigwedge^2V$, giving rise to a holomorphic curve
$f\wedge f'\:\Bbb C\to\Bbb P\bigl(\bigwedge^2V\spcheck\bigr)$.  We define
$\widebar T_{\bold y\wedge\bold y'}(r)$ and $N_{\bold y\wedge\bold y'}(r)$
analogously to $\widebar T_{\bold y}(r)$ and $N_{\bold y}(r)$.

Finally, let $\bold L$ be a finite set, all of whose elements are
$(n+1)$\snug-tuples $(L_0,\dots,L_n)$ of linearly independent linear forms
on $V$.  For each $z\in\Bbb C$ pick an element
$$(L_{z,0},\dots,L_{z,n})\in\bold L\;.$$
We also write $L_{r,\theta,i}=L_{re^{\sqrt{-1}\theta},i}$ for all $i$, let
$$\lambda_{r,\theta,i}(\bold v)
  = -\log\frac{|L_{r,\theta,i}(\bold v)|}{|\bold v|}\tag\05.1$$
be a Weil function for $L_{r,\theta,i}$, and let
$$m_{1,\bold y}(\bold L,r)
  = \frac1{n+1}\int_0^{2\pi}\sum_{i=0}^n
    \lambda_{r,\theta,i}(\bold y(re^{\sqrt{-1}\theta}))\frac{d\theta}{2\pi}\;.
  \tag\05.2$$
Of course this depends not only on $\bold L$ but also on the chosen function
$\Bbb C\to\bold L$.  Note also that this is independent of the lifting $\bold y$
of $f$, and that up to $O(1)$ it is independent of multiplying the $L_i$
(and therefore the $L_{r,\theta,i}$) by constants.

To help make sense of (\05.2), we note that if $\bold L$ has only one
element $(L_0,\dots,L_n)$, then
$$m_{1,\bold y}(\bold L,r) = \frac1{n+1}\sum_{i=0}^n m_f(\{L_i=0\},r)\;.$$
(Note, though, that having only one element in $\bold L$ corresponds to a
situation in which Theorem \00.2 follows trivially from the First Main Theorem.)

By (\05.1), we also have
$$m_{1,\bold y}(\bold L,r) - \widebar T_{\bold y}(r)
  = \frac1{n+1}\int_0^{2\pi}\sum_{i=0}^n
    -\log\bigl|L_{r,\theta,i}(\bold y(re^{\sqrt{-1}\theta}))\bigr|
      \frac{d\theta}{2\pi}\;.$$

If $I=\{i,j\}$ is a two-element subset of $\{0,\dots,n\}$, then we
let $L_{z,I}$ denote the linear form $L_{z,i}\wedge L_{z,j}$ on $\bigwedge^2V$,
let $\lambda_{z,I}$ be the corresponding Weil function, and let
$\lambda_{r,\theta,I}=\lambda_{re^{\sqrt{-1}\theta},I}$ as before.

\defn{\05.3}  We say that a collection $\Cal I$ of two-element subsets
of $\{0,\dots,n\}$ is {\bc balanced} if each $i\in\{0,\dots,n\}$ occurs
in the same number of elements of $\Cal I$.
\endit

If $\Cal I$ is a nonempty collection of two-element subsets of $\{0,\dots,n\}$,
then we define
$$m_{\Cal I,\bold y\wedge\bold y'}(\bold L,r)
  = \frac1{\#\Cal I}\int_0^{2\pi}\sum_{I\in\Cal I}
    \lambda_{r,\theta,I}
      (\bold y(re^{\sqrt{-1}\theta})\wedge\bold y'(re^{\sqrt{-1}\theta}))
    \frac{d\theta}{2\pi}\;,\tag\05.4$$
which as before reduces to
$$\frac1{\#\Cal I}\sum_{I\in\Cal I} m_{f\wedge f'}(\{L_I=0\},r)$$
in the case $\bold L=\{(L_0,\dots,L_n)\}$.

\lemma{\05.5}  Let $\Cal I$ be a nonempty balanced collection of two-element
subsets of $\{0,\dots,n\}$.  Then
$$2m_{1,\bold y}(\bold L,r) - m_{\Cal I,\bold y\wedge\bold y'}(\bold L,r)
  \le_\exc 2\,\widebar T_{\bold y}(r) - \widebar T_{\bold y\wedge\bold y'}(r)
    + O(\log^{+}T_f(r)) + o(\log r)\;.\tag\05.5.1$$
\endit

\demo{Proof}  We may assume that the coordinate functions $y_i$ of $\bold y$
never vanish simultaneously.  Indeed, if all coordinates of $\bold y$
are divisible by an entire function $g$, then all coordinates of
$\bold y\wedge\bold y'$ are divisible by $g^2$, so dividing $\bold y$ by $g$
leaves both sides of (\05.5.1) unchanged.

Since we now have $\widebar T_{\bold y}(r)=T_f(r)$, (\05.5.1) is equivalent
to
$$\split & 2m_{1,\bold y}(\bold L,r)
    - m_{\Cal I,\bold y\wedge\bold y'}(\bold L,r)
    + N_{\bold y\wedge\bold y'}(r) \\
  &\qquad \le_\exc 2\,T_f(r) - T_{f\wedge f'}(r)
    + O(\log^{+}T_f(r)) + o(\log r)\;.
  \endsplit\tag\05.5.2$$

The strategy of the proof is to apply McQuillan's Theorem \01.6
to $P=\Bbb P(V\spcheck)$.

First consider the first two height terms on the right-hand side.
The holomorphic map
$$\Bbb C @>(f,f\wedge f')>> \Fl_{1,2}(V) @>\phi>> \Bbb P(\Omega_{P/\Bbb C})$$
coincides with $f'$, where $\phi$ is the isomorphism of Corollary \03.8.
By the second assertion of Corollary \03.8, we then have
$$2\,T_f(r) - T_{f\wedge f'}(r) = -T_{\Cal O(1),f'}(r) + O(1)\;.\tag\05.5.3$$

We now consider the two proximity terms on the left-hand side of (\05.5.2).
Let $D_1,\dots,D_\ell$ be the divisors associated to the elements of $\bold L$,
and let $\mu_j$ be as in Theorem \01.6.  Then we claim that
$$2m_{1,\bold y}(\bold L,r) - m_{\Cal I,\bold y\wedge\bold y'}(\bold L,r)
  \le \int_0^{2\pi}\max_{1\le j\le\ell}
    \mu_j(f'(re^{\sqrt{-1}\theta}))\frac{d\theta}{2\pi} + O(1)\;,\tag\05.5.4$$
where the constant in $O(1)$ depends only on $\bold L$.  From the definitions,
we have
$$\split
  & 2m_{1,\bold y}(\bold L,r) - m_{\Cal I,\bold y\wedge\bold y'}(\bold L,r) \\
  &\qquad= \int_0^{2\pi} \biggl( \frac2{n+1}\sum_{i=0}^n
     \lambda_{r,\theta,i}(\bold y(re^{\sqrt{-1}\theta})) \\
    &\qquad\qquad- \frac1{\#\Cal I}\sum_{I\in\Cal I}
    \lambda_{r,\theta,I}
      (\bold y(re^{\sqrt{-1}\theta})\wedge\bold y'(re^{\sqrt{-1}\theta}))
    \frac{d\theta}{2\pi}\biggr)\;.\endsplit$$
Fix $r$ and $\theta$, let $L_i=L_{re^{\sqrt{-1}\theta},i}$ for all $i$,
let $\lambda_i=\lambda_{r,\theta,i}$ for all $i$, let
$\lambda_I=\lambda_{r,\theta,I}$ for all $I$, and choose $j$ so that $D_j$
is the divisor associated to $L_0,\dots,L_n$.  Then, to prove (\05.5.4),
it will suffice to show that
$$\split& \frac2{n+1}\sum_{i=0}^n \lambda_i(\bold y(re^{\sqrt{-1}\theta}))
    - \frac1{\#\Cal I}\sum_{I\in\Cal I}
    \lambda_I(\bold y(re^{\sqrt{-1}\theta})\wedge\bold y'(re^{\sqrt{-1}\theta}))
    \\
  &\qquad \le \mu_j(f'(re^{\sqrt{-1}\theta})) + O(1)\;,\endsplit\tag\05.5.5$$
with a constant depending only on $L_0,\dots,L_n$.  All of the terms in this
expression are unchanged up to $O(1)$ by a change of coordinates on $V$,
so we may assume that $L_i(X_0,\dots,X_n)=X_i$ for all $i$.  We also permute
indices such that
$$|y_0(re^{\sqrt{-1}\theta})|\ge|y_i(re^{\sqrt{-1}\theta})|
  \qquad\text{for all $i$}\;.$$
Then, letting $z_i=y_i(re^{\sqrt{-1}\theta})/y_0(re^{\sqrt{-1}\theta})$
and $z_i'=(y_i/y_0)'(re^{\sqrt{-1}\theta})$ for $i=0,\dots,n$,
we have $|z_i|\le1$ for all $i$,
$$\lambda_i(\bold y(re^{\sqrt{-1}\theta})) = -\log|z_i| + O(1)$$
for all $i$, and
$$\lambda_I(\bold y(re^{\sqrt{-1}\theta})\wedge\bold y'(re^{\sqrt{-1}\theta}))
    = -\log\frac{\vmatrix z_{i_1} & z_{i_2} \\ z_{i_1}' & z_{i_2}' \endvmatrix}
      {\max_{1\le k\le n}|z_k'|} + O(1)$$
for all $I=\{i_1,i_2\}$.  Thus the left-hand side of (\05.5.5) equals
$$\split& -\frac2{n+1}\sum_{i=0}^n \log|z_i|
    + \frac1{\#\Cal I}\sum_{\{i_1,i_2\}\in\Cal I}
      \log\frac{\vmatrix z_{i_1} & z_{i_2} \\ z_{i_1}' & z_{i_2}' \endvmatrix}
      {\max_{1\le k\le n}|z_k'|} + O(1) \\
  &\qquad= \frac1{\#\Cal I}\sum_{\{i_1,i_2\}\in\Cal I}
      \log\frac{\vmatrix z_{i_1} & z_{i_2} \\ z_{i_1}' & z_{i_2}' \endvmatrix}
      {|z_{i_1}||z_{i_2}|\max_{1\le k\le n}|z_k'|} + O(1)\;.\endsplit$$
since $\Cal I$ is balanced.  Substituting this and (\01.3.1) into (\05.5.5),
it follows that to prove (\05.5.5) it will suffice to show that
$$\log\frac{\vmatrix z_{i_1} & z_{i_2} \\ z_{i_1}' & z_{i_2}' \endvmatrix}
    {|z_{i_1}||z_{i_2}|\max_{1\le k\le n}|z_k'|}
  \le -\frac12\log\frac{|z_1'|^2+\dots+|z_n'|^2}
    {|z_1'/z_1|^2+\dots+|z_n'/z_n|^2} + O(1)$$
for all $0\le i_1<i_2\le n$, with a constant in $O(1)$ depending only on $n$.
After exponentiating and noting that
$\max|z_k'|^2\gg\ll|z_1'|^2+\dots+|z_n'|^2$, this is equivalent to
$$\left|\frac{z_{i_2}'}{z_{i_2}} - \frac{z_{i_1}'}{z_{i_1}}\right|
  \ll \sqrt{\fracwithdelims||{z_1'}{z_1}^2+\dots+\fracwithdelims||{z_n'}{z_n}^2}
  \;,$$
which is easy to see.  Thus (\05.5.5) is proved, so (\05.5.4) holds.

Finally, we note that
$$N_{\bold y\wedge\bold y'}(r) = N_{\text{Ram},f}(r)\;.\tag\05.5.6$$
Indeed, we may suppose without loss of generality that $y_0(z_0)\ne0$.
We have
$$\fracwithdelims(){y_i}{y_0}' = \frac{y_i'y_0-y_0'y_i}{y_0^2}
  = \frac{(\bold y\wedge\bold y')_{0i}}{y_0^2}\;.$$
Pick $i\in\{1,2,\dots,n\}$ such that $\ord_{z_0}(y_i/y_0)'$ is minimal.
This is the ramification order of $f$ at $z_0$.  It is also the order of
vanishing of the $0i$ coordinate of $\bold y\wedge\bold y'$; therefore
$N_{\bold y\wedge\bold y'}(r) \le N_{\text{Ram},f}(r)$.  The opposite
inequality is left to the reader (since it is not used here).

By (\05.5.3), (\05.5.4), (\05.5.6), and Theorem \01.6, we then have
$$\split& - 2\,T_f(r) + T_{f\wedge f'}(r)
  + 2m_{1,\bold y}(\bold L,r) - m_{\Cal I,\bold y\wedge\bold y'}(\bold L,r)
  + N_{\bold y\wedge\bold y'}(r) \\
  &\qquad\le T_{\Cal O(1),f'}(r)
  + \int_0^{2\pi}\max_{1\le j\le\ell}
    \mu_j(f'(re^{\sqrt{-1}\theta}))\frac{d\theta}{2\pi}
  + N_{\text{Ram},f}(r) + O(1) \\
  &\qquad\le_\exc O(\log^{+}T_f(r)) + o(\log r)\;.\endsplit$$
This gives (\05.5.2), as was to be shown.\qed
\enddemo

\beginsection{\06}{Main Step of the Proof}

We begin with some notation.  Let $V=\Bbb C^{n+1}$ as before, and let
$\bold x\:\Bbb C\to V$ be a holomorphic map whose coordinate functions
$x_0,\dots,x_n$ are linearly independent over $\Bbb C$.  This corresponds
to a holomorphic curve $f\:\Bbb C\to\Bbb P(V\spcheck)$ whose image is not
contained in any hyperplane.

As usual we let $\bold x^{(j)}=(x_0^{(j)},\dots,x_n^{(j)})$ be the
$j^{\text{th}}$ derivative of $\bold x$ ($j\in\Bbb N$), and following
Ahlfors we let
$$\align X^d &= \bold x \wedge \bold x' \wedge \dots \wedge \bold x^{(d-1)}
    \in \bigwedge\nolimits^d V\;, \\
  \widebar T_{d,\bold x}(r)
    &= \int_0^{2\pi}\log\bigl|X^d(re^{\sqrt{-1}\theta})\bigr|
      \frac{d\theta}{2\pi}\;,\endalign$$
and let $N_{d,\bold x}(r)$ be the counting function for the simultaneous
vanishing of the coordinates of $X^d$, for $d=0,\dots,n+1$.  (If $d=0$ then
$X^0\:\Bbb C\to\bigwedge^0 V=\Bbb C$ is the constant map $1$.)  Let $F^d$
denote the corresponding map to $\Bbb P(\bigwedge^d V\spcheck)$; we then have
$$T_{\Cal O(1),F^d}(r) = \widebar T_{d,\bold x}(r) - N_{d,\bold x}(r)\;.$$
We also write $T_{d,f}(r)=T_{\Cal O(1),F^d}(r)$.  Note that $T_{1,f}(r)=T_f(r)$
and $T_{0,f}(r)=0$ for all $r$.

Let $\bold L$, $L_{z,i}$, and $L_{r,\theta,i}$ be as before.  If $I$
is a $d$\snug-element subset of $\{0,\dots,n\}$, then we let
$L_{z,I}=L_{z,i_1}\wedge\dots\wedge L_{z,i_d}$, where $I=\{i_1,\dots,i_d\}$
with $i_1<\dots<i_d$.  Also let $\lambda_{z,I}$ be the corresponding Weil
function on $\Bbb P(\bigwedge^dV\spcheck)$, and let
$\lambda_{r,\theta,I}=\lambda_{re^{\sqrt{-1}\theta},I}$ as before.
Extending (\05.2) in a manner similar to (\05.4), we let
$$m_{d,f}(\bold L,r) = \binom{n+1}{d}^{-1}\int_0^{2\pi}\sum_{\#I=d}
  \lambda_{r,\theta,I}(F^d(re^{\sqrt{-1}\theta}))\frac{d\theta}{2\pi}\;.
  \tag\06.1$$

The main step in the proof is then as follows.

\prop{\06.2}  For all $d=1,\dots,n$, we have
$$\split& -m_{d-1,f}(\bold L,r) + 2m_{d,f}(\bold L,r) - m_{d+1,f}(\bold L,r) \\
  &\qquad\le_\exc -\widebar T_{d-1,\bold x}(r) + 2\widebar T_{d,\bold x}(r)
    - \widebar T_{d+1,\bold x}(r)
    + O(\log^{+}T_{d,f}(r)) + o(\log r)\;.\endsplit\tag\06.2.1$$
\endit

\demo{Proof}  The proof works by applying Lemma \05.5 to $\bold y:=X^d$.

We start by noting that the derivative of
$\bold y = X^d = \bold x\wedge\bold x'\wedge\dots\wedge\bold x^{(d-1)}$
can be computed by a Leibniz-like relation.  This gives $d$ terms, all
but the last of which vanish, giving
$$\bold y' = \bold x \wedge \bold x' \wedge \dots \wedge \bold x^{(d-2)}
  \wedge \bold x^{(d)}\;.$$

When applying Lemma \05.5, we will use linear forms on $\bigwedge^dV$
obtained as wedge products of $d$ linear forms on $V$.  When working with
these forms, the classical formula
$$(L_1\wedge\dots\wedge L_d)(\bold x_1\wedge\dots\wedge\bold x_d)
  = \det(L_i(\bold x_j))_{1\le i,j\le n}\tag\06.2.2$$
is useful.

We will also use the following formula, also used by Ahlfors and Schmidt;
see \cite{V~1, Lemma~6.3.14}.  Let $A$ be a $(d-1)\times(d-1)$ matrix,
let $B$ and $B'$ be $(d-1)\times1$ matrices, let $C$ and $C'$ be $1\times(d-1)$
matrices, and let $d$, $e$, $f$, and $g$ be scalars.  Then
$$\vmatrix\vmatrix A & B \\ C & d \endvmatrix
    & \vmatrix A & B' \\ C & e \endvmatrix \\
    \vspace{1.5\jot}
    \vmatrix A & B \\ C' & f \endvmatrix
    & \vmatrix A & B' \\ C' & g \endvmatrix \endvmatrix
  = \vmatrix A \endvmatrix
    \vmatrix A & B & B' \\ C & d & e \\ C' & f & g \endvmatrix\;.$$
Combining this with (\06.2.2) gives
$$(L_{z,I}\wedge L_{z,J})(\bold y\wedge\bold y')
  = \vmatrix L_{z,I}(\bold y) & L_{z,J}(\bold y) \\
    L_{z,I}(\bold y') & L_{z,J}(\bold y') \endvmatrix
  = \pm L_{z,I\cap J}(X^{d-1})L_{z,I\cup J}(X^{d+1})\;,\tag\06.2.3$$
where $I$ and $J$ are $d$\snug-element subsets of $\{0,\dots,n\}$ that
{\it differ by exactly one element.}

Therefore, for $I,J\subseteq\{0,\dots,n\}$ with $\#I=\#J=d$, we let
$$\dist(I,J)=\#(I\setminus(I\cap J))=\#(J\setminus(I\cap J))$$
be the number of elements in which they differ.  Let $\Cal I$ be the
collection of sets $\{I,J\}$ with $I,J\subseteq\{0,\dots,n\}$, $\#I=\#J=d$,
and $\dist(I,J)=1$.  Then, with a suitable definition of $\bold L'$, we have
$$m_{\Cal I,\bold y\wedge\bold y'}(\bold L',r)
  = \frac1{\#\Cal I}\int_0^{2\pi}\sum\Sb\#I=\#J=d\\\dist(I,J)=1\endSb
    \lambda_{r,\theta,\{I,J\}}
      (\bold y(re^{\sqrt{-1}\theta})\wedge\bold y'(re^{\sqrt{-1}\theta}))
    \frac{d\theta}{2\pi}$$
and therefore, by (\05.1) and (\06.2.3),
$$\split& \widebar T_{\bold y\wedge\bold y'}(r)
    - m_{\Cal I,\bold y\wedge\bold y'}(\bold L',r) \\
  &\qquad= \frac1{\#\Cal I}\int_0^{2\pi}\sum\Sb\#I=\#J=d\\\dist(I,J)=1\endSb
    \log\bigl|(L_{r,\theta,I}\wedge L_{r,\theta,J})
      (\bold y(re^{\sqrt{-1}\theta})\wedge\bold y'(re^{\sqrt{-1}\theta}))\bigr|
    \frac{d\theta}{2\pi} \\
  \allowdisplaybreak
  &\qquad= \frac1{\#\Cal I}\int_0^{2\pi}\sum\Sb\#I=\#J=d\\\dist(I,J)=1\endSb
    \Bigl(\log\bigl|L_{r,\theta,I\cap J}(X^{d-1}(re^{\sqrt{-1}\theta}))\bigr| \\
      &\qquad\qquad+ \log\bigl|L_{r,\theta,I\cup J}
        (X^{d+1}(re^{\sqrt{-1}\theta}))\bigr|\Bigr)\frac{d\theta}{2\pi} \\
  \allowdisplaybreak
  &\qquad= \binom{n+1}{d-1}^{-1}\int_0^{2\pi}\sum_{\#I=d-1}
    \log\bigl|L_{r,\theta,I}(X^{d-1}(re^{\sqrt{-1}\theta}))\bigr|
      \frac{d\theta}{2\pi} \\
    &\qquad\qquad+ \binom{n+1}{d+1}^{-1}\int_0^{2\pi}\sum_{\#I=d+1}
    \log\bigl|L_{r,\theta,I}(X^{d+1}(re^{\sqrt{-1}\theta}))\bigr|
      \frac{d\theta}{2\pi} \\
  &\qquad= \widebar T_{d-1,\bold x}(r) - m_{d-1,f}(\bold L,r)
    + \widebar T_{d+1,\bold x}(r) - m_{d+1,f}(\bold L,r)\;.
  \endsplit$$
Here we also used the fact that, as $\{I,J\}$ varies over $\Cal I$,
the intersection $I\cap J$ varies over all $(d-1)$\snug-element subsets of
$\{0,\dots,n\}$ with equal frequency, and the same holds for the union
$I\cup J$.

Straight from the definitions we also have
$$m_{1,\bold y}(\bold L',r) = m_{d,f}(\bold L,r)\;,
  \qquad \widebar T_{\bold y}(r) = \widebar T_{d,\bold x}(r)\;,
  \qquad\text{and}\qquad T_{\bold y}(r) = T_{d,f}(r)\;,$$
and therefore (\06.2.1) follows from (\05.5.1).\qed
\enddemo

\prop{\06.3}  In (\06.2.1), the error term $O(\log^{+}T_{d,f}(r))$ can be
replaced by\break
$O(\log^{+}T_f(r))$.
\endit

\demo{Proof}  Let $1\le d\le n$ and let $G=\Gr_d(V)$.  The composite map
$$\Bbb C @>(X^{d-1},X^d,X^{d+1})>> \Fl_{d-1,d,d+1}(V)
  @>\phi>> \Bbb P(\Omega_{G/\Bbb C})\;,$$
where $\phi$ is the map of Proposition \03.7, is just the map $(F^d)'$.
Applying McQuillan's Theorem \01.2 with $D=0$ gives
$$T_{\Cal O(1),(F^d)'}(r) \le_\exc O(\log^{+}T_{d,f}(r)) + o(\log r)\;.$$
By Proposition \03.7, $\phi^{*}\Cal O(1)\cong\Cal O(1,-2,1)$, so by
functoriality of the height (characteristic) function,
$$T_{d-1,f}(r) - 2T_{d,f}(r) + T_{d+1,f}(r)
  \le_\exc O(\log^{+}T_{d,f}(r)) + o(\log r)\;.$$
By induction on $d$, we then have
$$T_{d,f}(r) \le_\exc 2^{d-1}T_f(r) + O(\log^{+}T_f(r)) + o(\log r)\;.
  \tag\06.3.1$$
This implies the result.\qed
\enddemo

\beginsection{\07}{Conclusion of the Proof}

Let $f\:\Bbb C\to\Bbb P^n$ be a holomorphic curve whose image is not contained
in any hyperplane.  Lift it to a holomorphic map
$\bold x\:\Bbb C\to\Bbb C^{n+1}=V$ such that the coordinate functions
$x_0,\dots,x_n$ have no common zeroes.  Then we have $N_{1,f}(r)=0$
for all $r$, so
$$\widebar T_{1,\bold x}(r) = T_f(r)\;.$$
Also, since $\bigwedge^{n+1}V\cong\Bbb C$, we have
$$T_{n+1,f}(r) = O(1)\;.$$

Now let $H_1,\dots,H_q$ be the hyperplanes in Theorem \00.2.  Adding
hyperplanes only strengthens the inequality, so we may assume that
$\bigcap H_i=\emptyset$.  The Weil functions $\lambda_{H_j}$ are bounded
from below, so the left-hand side of (\00.2.1) is changed by only $O(1)$
if we assume that all subsets $J$ have at least $n+1$ elements.  Also,
given a subset of $H_1,\dots,H_q$ in general position, no point in $\Bbb P^n$
is close to more than $n$ hyperplanes in the subset.  Therefore, again up
to $O(1)$, we may assume that all subsets $J$ have exactly $n+1$ elements.

Now choose linear forms $L_1,\dots,L_q$ defining the hyperplanes
$H_1,\dots,H_q$, respectively, and let $\bold L$ be the collection of
all tuples $(L_{i_0},\dots,L_{i_n})$
such that $1\le i_0<i_1<\allowmathbreak \dots<i_n\le q$ and
$H_{i_0},\dots,H_{i_n}$ are in general position.  For each $z\in\Bbb C$ pick
$(L_{z,0},\dots,L_{z,n})\in\bold L$ such that
$$\max_J\sum_{j\in J} \lambda_{H_j}(f(z))
  = \sum_{i=0}^n \lambda_{H_{z,i}}(f(z))\;,$$
where $H_{z,i}$ is the hyperplane determined by $L_{z,i}$.  We then have
$$\int_0^{2\pi}\max_J\sum_{i\in J}\lambda_{H_i}(f(re^{\sqrt{-1}\theta}))
    \frac{d\theta}{2\pi}
  = (n+1)m_{1,f}(\bold L,r) + O(1)\;,\tag\07.1$$
where the sets $J$ are as in (\00.2.1).

By Propositions \06.2 and \06.3, we have
$$\split (n+1)m_{1,f}(\bold L,r)
  &= -nm_{0,f}(\bold L,r) + (n+1)m_{1,f}(\bold L,r) - m_{n+1,f}(\bold L,r) \\
  &= \sum_{d=1}^n (n+1-d)\bigl(-m_{d-1,f}(\bold L,r) + 2m_{d,f}(\bold L,r)
    - m_{d+1,f}(\bold L,r)\bigr) \\
  &\le_\exc \sum_{d=1}^n (n+1-d)\bigl(-\widebar T_{d-1,f}(r)
    + 2\widebar T_{d,f}(r) - \widebar T_{d+1,f}(r)\bigr) \\
    &\qquad+ O(\log^{+}T_f(r)) + o(\log r) \\
  &= -n\,\widebar T_{0,f}(r) + (n+1)\widebar T_{1,f}(r) - \widebar T_{n+1,f}(r)
    \\
    &\qquad+ O(\log^{+}T_f(r)) + o(\log r) \\
  &= (n+1)T_f(r) - N_{n+1,f}(r) + O(\log^{+}T_f(r)) + o(\log r)\;.
  \endsplit\tag\07.2$$
Combining (\07.1) and (\07.2) gives (\00.2.1), upon noting that
$X^{n+1}\:\Bbb C\to\bigwedge^{n+1}V=\Bbb C$ is exactly the Wronskian used
in (\00.2.1).\qed

\enddocument